\documentstyle[12pt]{article}
\date{}
\newtheorem{th}{Theorem}[section]
\newtheorem{prop}{Proposition}[section]
\newtheorem{cor}{Corollary}[section]
\newtheorem{lem}{Lemma}[section]

\makeatletter
\renewcommand{\@begintheorem}[2]{\par\it {\bf #1\ #2. }}
\renewcommand{\@endtheorem}{}

\makeatother

\title{Hilbert modules over locally $C^*$-algebras}
\author{Yu. I. Zhuraev, F. Sharipov
\thanks{ The work is partially supported by the INTAS grant 96-1099. }}

\begin{document}
\maketitle

\begin{abstract}
In the present paper the notion of a Hilbert module over a locally
$C^*$-algebra is discussed and some results are obtained on this
matter. In particular, we give a detailed proof of the known result that
the set of adjointable endomorphisms of such modules is itself a locally
$C^*$-algebra.

\end{abstract}

 \begin{flushleft}
\Large\bf Introduction
 \end{flushleft}

General theory of Hilbert modules over an arbitrary $C^*$-algebra
was constructed in the papers of
W.~Paschke~\cite{1} and M.~Rieffel~\cite{2}
as a natural generalization of the Hilbert spaces theory.
This generalization arises if a $C^*$-algebra takes place of
the field of scalars ${\bf C}$.
Theory of Hilbert modules is applied
to many fields of mathematics, in particular, to theory of
$C^*$-algebras \cite{3,24}, to theory of vector bundles \cite{12}, to
theory of index of elliptic operators \cite{4,11}, to $K$-theory
\cite{5,6}, to $KK$-theory of G.~Kasparov \cite{7}, to theory of quantum
groups and unbounded operators \cite{13}, to some physical problems
\cite{9,14} etc. There is also a number of papers dedicated to theory of
Hilbert $C^*$-modules proper (see, for example, \cite{10,15,16,17}).

Other topological $*$-algebras, for example, locally $C^*$-algebras and
some group algebras can be met in applications together with
$C^*$-algebras. By analogy with $C^*$-algebras locally $C^*$-algebras
are applied to relativistic quantum mechanics~\cite{18,19}.
Therefore it seems useful to develop theory of Hilbert modules
over other topological $*$-algebras as well. As we know the
paper~\cite{20} of
A.~Mallios was the first one in this direction. In that paper finitely
generated modules equipped with inner products
over some topological $*$-algebras were considered in connection with
Hermitian $K$-theory, the standard Hilbert module $l_2(A)$ over a locally
$C^*$-algebra $A$ was introduced and the index theory for elliptic
operators over a locally $C^*$-algebra was constructed.

In the present paper the notion of a Hilbert module
over a locally $C^*$-algebra is discussed and some results
are obtained on this matter. In particular, it is proved that the set of
adjointable endomorphisms of such modules is itself a locally
$C^*$-algebra.

After the first version of this paper appeared we were informed that our
main result was already known to specialists \cite{a,a1,b}, but we think
that our detailed proof may still be of interest.

 \section{Locally $C^*$-algebras}
We start with some information from \cite{21,22,23} on locally
$C^*$-algebras.

{\bf Definition 1.1.}
 A complex algebra $A$ is called a {\em $LMC$-algebra\/} if
 it is a separable locally convex space with respect to some
 family of seminorms $\{P_\alpha\}_{\alpha \in \Delta}$
 satisfying the following condition:

(1) $P_\alpha(ab) \le P_\alpha(a)P_\alpha(b)$ for all $a,b \in A.$

An involutive $LMC$-algebra $A$ is called a {\em $LMC$ $*$-algebra\/} if
the following condition holds:

(2) $P_\alpha(a^*)=P_\alpha(a)$ for all  $a \in A$ and $\alpha\in\Delta.$

A complete $LMC$ $*$-algebra $A$ is called a {\em locally $C^*$-algebra\/}
({\em $LC^*$-algebra\/}) if

(3) $P_\alpha(a^*a)=(P_\alpha(a))^2$ for all  $a \in A$ and $\alpha \in
\Delta.$

A {\em family of $C^*$-seminorms\/} is the family of seminorms
$\{P_\alpha\}_{\alpha\in\Delta}$ satisfying condition (3)
of Definition 1.1.

Let us give some examples of locally $C^*$-algebras \cite{21,22,23}.

{\bf Example 1.1.}  Any $C^*$-algebra is a $LC^*$-algebra.

{\bf Example 1.2.} A closed $*$-subalgebra of a $LC^*$-algebra is
a $LC^*$-algebra.

{\bf Example 1.3.} Let $M$ be a completely regular $k$-space
(\cite{25}, p. 300) and $C(M)$ an algebra of all continuous complex-valued
functions on $M$. For each compact space $K \subset M$ put
$$
q_K(f):=\sup\limits_{x \in K} |f(x)|,\quad f \in C(M).
$$
Then the function $q_K$ is a $C^*$-seminorm on $C(M)$ and
with respect to the family of these seminorms $C(M)$ is a $LC^*$-algebra.

{\bf Example 1.4.} Let $\Lambda$ be a directed set of indices,
$\{H_\lambda\}_{\lambda\in\Lambda}$ a family of Hilbert spaces
such that $H_\lambda \subseteq H_\mu$  and
 $$
(\cdot,\cdot)_\lambda=(\cdot,\cdot)_\mu|_{H_\lambda}
$$
for $\lambda \le \mu$. Here $(\cdot,\cdot)_\lambda$ is an inner product
on $H_\lambda,\;\lambda \in \Lambda$.
Consider a locally convex space
$$
H:=\lim_{\longrightarrow}{}^{\phantom{\chi}}_{\lambda}\, H_\lambda=
\bigcup_\lambda H_\lambda.
$$
The space $H$ equipped with the topology of the inductive limit
is called {\em a locally Hilbert space.\/} We denote by $L(H)$
the set of all linear continuous operators $T:H \to H$ such that
$$T=\lim\limits_{\longrightarrow}{}^{\phantom{\chi}}_{\lambda}\,
T_\lambda,\quad T_\lambda \in B(H_\lambda).$$
Here $B(H_\lambda)$ is the space of all linear bounded operators
on $H_\lambda$ and $(T_\lambda)_{\lambda\in\Lambda}$ is the inductive
family of the operators $T_\lambda \in B(H_\lambda)$.
It is clear that $L(H)$ is an algebra. Furthemore, if
$(T_\lambda)_{\lambda \in \Lambda}$ is an inductive family of
linear bounded operators on $H_\lambda, \lambda \in\Lambda$,
then the family of adjoint operators $(T^*_\lambda)_{\lambda \in\Lambda}$
is inductive too. The map
$$
*:L(H) \to L(H),\quad T \mapsto T^*=
\lim_{\longrightarrow}{}^{\phantom{\chi}}_{\lambda}\, T^*_\lambda
$$
defines an involution on $L(H)$. If $||\cdot||_\lambda$ is the operator
norm on $B(H_\lambda)$ then the function
$$
q_\lambda(T):=||T_\lambda||_\lambda,\; T \in L(H)
$$
is a $C^*$-seminorm on $L(H)$ for each $\lambda \in \Lambda$
and $L(H)$ is a $LC^*$-algebra with respect to the family of seminorms
$\{q_\lambda\}_{\lambda \in \Lambda}$.

 \begin{th} {\rm (\cite{21}, Theorem 5.1).}
 Any $LC^*$-algebra is isomorphic to a closed $*$-subalgebra
 of $L(H)$ for some locally Hilbert space $H$.
 \end{th}

Let $A$ be a $LC^*$-algebra with respect to a family of
$C^*$-seminorms $\{P_\alpha\}_{\alpha \in \Delta}$.
We denote by $I_\alpha$ the kernel of the seminorm $P_\alpha$, i.e.
the set of elements $a \in A$ such that $P_\alpha(a)=0$.
It is clear that $I_\alpha$ is a closed $*$-ideal in $A$.
Therefore the quotient space $A_\alpha=A/I_\alpha$ is a normed
$*$-algebra with respect to the norm
$$
||a_\alpha||:=P_\alpha(a), \quad a_\alpha=a+I_\alpha \in A_\alpha.
$$
It follows from Theorem 2.4 of \cite{22} that the algebra $A_\alpha$
is complete, i.e. it is a $C^*$-algebra.

By $e$ we will denote the identity element of $A$. Clearly
$P_\alpha(e)=1$ for any non-zero seminorm $P_\alpha$.

If $A$ is an algebra without the identity element, then by $A^+$ we denote
its unitalization. By Theorem 2.3 of \cite{21},
any seminorm $P_{\alpha}$ can be extended up to a $C^*$-seminorm
$P_{\alpha}^+$ on $A^+$ and $A^+$ is a locally $C^*$-algebra
with respect to the family of seminorms $P_{\alpha}^+$.

The {\em spectrum\/} of an element $a$ of a unital $LC^*$-algebra $A$
is the set $Sp(a)=Sp_A(a)$ of complex numbers $z$ such that
$a-z\cdot1$ is not invertible. If $A$ is a non-unital algebra, then
the {\em spectrum\/} of an element $a \in A$ is its spectrum in the
$LC^*$-algebra $A^+$. It follows from Corollary 2.1 of \cite{21} that
  $$
  Sp_A(a)=\bigcup_{\alpha\in\Delta}Sp_{A_{\alpha}}(a_{\alpha}), \quad
  a_{\alpha}=a+I_{\alpha}\eqno(1.1)
  $$
  for each $a \in A$. An element $a \in A$ is called {\em positive\/}
 (and is written $a \ge 0$) if it is Hermitian, i.e. $a=a^*$
  and one of the following equivalent (see \cite{21}, Proposition 2.1)
conditions is true:

(1) $Sp(a) \subset [0,\infty);$

(2) $a=b^*b$ for some $b \in A;$

(3) $a=h^2$ for some Hermitian $h \in A.$

  Besides, the set of positive elements $P^{+}(A)$ is a closed
  convex cone in $A$ and $P^{+}(A) \bigcap (-P^{+}(A))=\{0\}.$

{\bf Remark 1.1.} If $a \in A$ is a positive element then
  there exists a unique positive element $h \in A$ satisfying
  the condition (3). This element is called a {\em square root\/}
  of $a$ and is denoted by $a^{\frac{1}{2}}=\sqrt{a}$.
  For elements $a,b \in A$ the inequality $a \ge b$ (or $b \le a$)
  means that $a-b \ge 0$.

   \begin{lem} {\rm (\cite{21}).}

(a) If $a,b \in P^+(A)$ and $a \le b$, then $P_{\alpha}(a)
\le P_{\alpha}(b)$ for all $\alpha \in \Delta.$

(b) If $e \in A, b \in A$ and $b \ge e$ then the element $b$ is
  invertible and $b^{-1} \le e.$

(c) If elements $a,b \in A$ are invertible and $0 \le a \le b$ then
  $a^{-1} \ge b^{-1}.$

(d) If $a,b,c \in A$ and $a \le b$  then $c^*ac \le c^*bc.$
   \end{lem}

   \begin{lem}
 Let $A$ be a unital locally $C^*$-algebra, $a \in P^+(A)$ and $t$ be a
positive number. Then for any $\alpha \in \Delta$ the following relations
  hold:

(a) $P_{\alpha}((e+ta)^{-1}) \le 1;$

(b) $P_{\alpha}(a(e+a)^{-1}) \le 1;$

(c) $P_{\alpha}(e-a) \le 1$ if $P_{\alpha}(a) \le 1.$
   \end{lem}

{\bf Proof.} Since $e+ta \ge e,$ we obtain from (b) and (a) of Lemma 1.1
 that the element $e+ta$ is invertible and
 $$
 P_{\alpha}((e+ta)^{-1}) \le  P_{\alpha}(e)=1.
 $$
 Using $e+a \ge a$ and statement (d) of Lemma 1.1 for
$c=\sqrt{(e+a)^{-1}}$ we obtain
 $$
 a(e+a)^{-1} \le e.
 $$
 This yields (b).

 We will prove (c). Since $A_{\alpha}=A/I_{\alpha}$ is a $C^*$-algebra
 with the identity element $e_{\alpha}=e+I_{\alpha}$ and
 $||a_{\alpha}||=P_{\alpha}(a) \le 1,$ we obtain that
 $$
 P_{\alpha}(e-a)=||e-a+I_{\alpha}||=||e_{\alpha}-a_{\alpha}|| \le 1.
 $$
 The lemma is proved.

We denote by $A^s$ the set of all elements $a \in A$ such that
 $$
 ||a||^s:=\sup_{\lambda \in \Delta}P_\alpha(a)<\infty.
 $$
 Then $A^s$ is a $*$-subalgebra of $A$ and $||\cdot||^s$ is a norm
  for $A^s$. Moreover, the following theorem is true.

  \begin{th} {\rm (\cite{22}, Theorem 2.3).}
 The algebra $A^s$ is dense in $A$ and is a $C^*$-algebra with
 respect to the norm $||\cdot||^s$.
  \end{th}

 An {\em approximate identity\/} of a locally $C^*$-algebra $A$ is
 any increasing net $\{u_{\lambda}\}_{\lambda \in \Lambda}$
 of positive elements such that

1) $P_{\alpha}(u_{\lambda}) \le 1$ for all $\alpha \in \Delta, \lambda \in
\Lambda;$

2) $\lim \limits_{\lambda}(a-au_{\lambda})=\lim
\limits_{\lambda}(a-u_{\lambda}a)=0$ for any $a \in A.$

Any locally $C^*$-algebra (and its closed ideal) has an
  approximate identity (see \cite{21}, Theorem 2.6).

  \section{Hilbert modules over $LC^*$-algebras}
Let $A$ be a $LC^*$-algebra with respect to the family of
 $C^*$-seminorms $\{P_\alpha \}_{\alpha \in \bigtriangleup}.$
 Below we will assume that the algebra $A$ has the unit $e$.
 Considering $A^+$ instead of $A$, one can easily extend all results
 for the case of non-unital algebras.

{\bf Definition 2.1.} Let $X$ be a right $A$-module.
An {\em $A$-valued inner product\/} on $X$ is a map
$<.,.>: X\times X \to A$ satisfying the following conditions:

$(1) <x, x> \ge 0$ for all $x \in X$;

$(2) <x,x>=0$ if and only if $x=0$;

$(3) <x_1+x_2,  y>=<x_1, y>+<x_2, y>$ for all $x_1, x_2, y \in X$;

$(4) <x a, y>=<x,y>a$ for all $x, y \in X, a \in A$;

$(5) <x,y>^*= <y, x>$ for all $x, y \in X$.

A right $A$-module equipped with an $A$-valued inner product
is called a {\em pre-Hilbert $A$-module.\/}

 \begin{lem}
  Let $X$ be a pre-Hilbert $A$-module. Then for each
$\alpha \in \bigtriangleup$ and for all $x,y \in X$
the following Cauchy-Bunyakovskii inequality holds:
$$ P_\alpha(<x,y>)^2 \le P_\alpha(<x,x>) P_\alpha(<y,y>).
\eqno(2.1) $$
 \end{lem}

{\bf Proof.} Suppose $x,y \in X$, $b \in A.$ Let us consider the
expression
$$ <x+yb, x+y b>=<x,x>+b^*<x,y>+<y,x>b
+b^*<y,y>b \ge 0. \eqno(2.2) $$
Assuming $P_\alpha(<y,y>) \ne 0$, we put
$$
b=-\frac{<x,y>}{P_\alpha(<y,y>)}
$$
in (2.2). It now follows that
$$<x,x>-
\frac{<y,x><x,y>}{P_{\alpha}(<y,y>)}-\frac{<y,x><x,y>}{P_{\alpha}(<y,y>)}+
\frac{<y,x><y,y><x,y>}{P_{\alpha}(<y,y>)^2} \ge 0,$$
whence, using item (a) of Lemma 1.1, we obtain
$$P_{\alpha}\left(\frac{2<y,x><x,y>}{P_{\alpha}(<y,y>)}\right) \le
P_{\alpha}\left(<x,x>+\frac{<y,x><y,y><x,y>}{P_{\alpha}(<y,y>)^2}\right).$$
Therefore,
$$
\frac{2P_{\alpha}(<x,y>)^2}{P_{\alpha}(<y,y>)}  \le
P_{\alpha}(<x,x>)+\frac{P_{\alpha}(<y,x>)P_{\alpha}(<y,y>)
P_{\alpha}(<x,y>)}{P_\alpha(<y,y>)^2}
$$
$$
=P_{\alpha}(<x,x>)+\frac{P_\alpha(<x,y>)^2}{P_\alpha(<y,y>)}.
$$
This implies inequality (2.1).
If $P_\alpha(<x,x>) \ne 0$, then it is true by the same reason.

Let now $P_\alpha(<x,x>)= P_\alpha(<y,y>)=0$.
Putting $b=-<x,y>$ in $(2.2)$, we get
$$
<x,x>+<y,x><y,y><x,y> \ge 2<y,x><x,y>,
$$
whence
$$
2P_\alpha(<x,y>)^2=2P_\alpha(<y,x><x,y>) \le P_\alpha(<x,x>)+
$$
$$
+P_\alpha(<y,x>) P_\alpha(<y,y>) P_\alpha(<x,y>)=0.
$$
Thus $P_\alpha(<x,y>)=0$ and inequality (2.1) is true in this case too.
The lemma is proved.


 \begin{lem}\label{Lemma2.2}
 Let $X$ be a pre-Hilbert $A$-module with respect to an inner product
$<.,.>$.  Put $$
\bar P_\alpha(x):= P_\alpha (<x,x>)^{\frac12}      \eqno(2.3)
$$
for any $\alpha \in \Delta$. Then the function $\bar P_\alpha$
is a seminorm on $X$ and the following conditions hold:

(1) $\bar P_\alpha(x a) \le \bar P_\alpha(x) P_\alpha(a)$ for all $x \in
X,\; a \in A;$

(2) if $\bar P_\alpha(x)= 0$ for all $\alpha \in \Delta$, then $x= 0;$

(3) $\bar P_\alpha(x)=\sup_{\bar P_\alpha(y) \le 1}P_\alpha(<x,y>)$ for
all $x \in X, \; \alpha \in \Delta.$
 \end{lem}

 {\bf Proof.} Under the axioms of seminorm we have
$$ \bar P_\alpha(x)= P_\alpha (<x,x>)^\frac12 \ge 0$$
and
$$ \bar P_\alpha(\lambda
 x)=\sqrt{P_\alpha(<\lambda x,\lambda x>)}= \sqrt{P_\alpha(\bar \lambda
 <x,
x>\lambda)} $$ $$ =\sqrt{P_\alpha(|\lambda|^2<x,x>)}=\sqrt{|\lambda|^2
P_\alpha(<x,x>)}= |\lambda|\bar P_\alpha (x)$$
for all $x \in X$ and $\lambda \in {\bf C}$.

Using the Cauchy-Bunyakovskii inequality one has
 $$
\bar P_\alpha (x+y)=\sqrt{P_\alpha(<x+y, x+y>)}
$$
$$
=\sqrt{P_\alpha(<x,x>+<y,x>+<x,y>+<y,y>)}
$$
$$
\le \sqrt{P_\alpha(<x,x>)+2 P_\alpha(<x,y>)+P_\alpha(<y,y>)}
$$
$$
\le \sqrt{\bar P_\alpha (x)^2+2 \bar P_\alpha (x) \bar P_\alpha (y)+
\bar P_\alpha (y)^2}= \bar P_\alpha (x)+\bar P_\alpha (y)
$$
for all $x, y \in X$.
Thus $\bar P_\alpha$ is a seminorm on $X$.

Let us prove (1). For all $a \in A$ and $x \in X$ we derive
$$
\bar P_\alpha(xa)=\sqrt{ P_\alpha(<xa,xa>)}=
\sqrt{ P_\alpha(a^* <x,x> a)}
\le \sqrt{ P_\alpha(a^*) P_\alpha(<x,x>) P_\alpha(a)}= P_\alpha(a)
\bar P_\alpha(x).
$$

Suppose $x \in X$ and $\bar P_\alpha(x)=0$ for all
$\alpha \in \bigtriangleup$. Then $P_\alpha(<x,x>)=0$ for all
 $\alpha \in \bigtriangleup$. Therefore $<x,x>=0$ and, consequently,
 $x=0$. Thus (2) is true. Equality (3) follows easy from
 the Cauchy-Bunyakovskii inequality. The lemma is proved.

 Lemma \ref{Lemma2.2} implies that $X$ is a separable locally convex space
 with respect to the family of seminorms
$\{\bar P_\alpha:\alpha \in \bigtriangleup\}$.

{\bf Definition 2.2.} Let $X$ be a pre-Hilbert $A$-module equipped with
 the inner product $<.,.>$. If $X$ is a complete locally convex space
 with respect to the family of seminorms
$\{\bar P_\alpha\}_{\alpha \in \bigtriangleup}$ defined by (2.3),
 then it is called a {\em Hilbert $A$-module.\/}

{\bf Example  2.1.} Any closed right ideal $I$ of a locally $C^*$-algebra
$A$ equipped with the inner product $<a,b>=a^* b$ is a Hilbert $A$-module.

{\bf Example 2.2} (\cite{20}).
Let $l_2 (A)$ be the set of all sequences
$x=(x_n)_{n \in {\bf N}}$ of elements from a locally $C^*$-algebra $A$
 such that the series
$$
\sum\limits_{i=1}^{\infty} x^*_i x_i
$$
is convergent in $A$. Then $l_2 (A)$ is a right Hilbert $A$-module
with respect to the pointwise operations and the inner product
$$
<x,y>=\sum\limits_{i=1}^{\infty} x^*_i y_i.
$$

Let $A$ be a locally $C^*$-algebra and $X$ a right Hilbert $A$-module
equipped with the inner product $<.,.>$. We will denote by $X^s$ the set
of all $x \in X$ such that $<x,x> \in A^s$. It is verified immediately
that $X^s$ is $A^s$-module.

 \begin{lem}
For all $x \in X$ and $\alpha \in \Delta$
one has:

a) $\bar P_\alpha(x(e+\sqrt{<x,x>})^{-1}) \le 1;$

b) $\lim \limits_{t \to +0} \bar P_\alpha(x-x(e+ t\sqrt{<x,x>})^{-1})=0.$
 \end{lem}

{\bf Proof.} a) Let $x \in X$ be an arbitrary element. Then for each
$\alpha \in \Delta$  we have by statement (b) of Lemma 1.2 that
$$
\bar P_\alpha(x(e+\sqrt{<x,x>})^{-1})^2=
P_\alpha(<x(e+\sqrt{<x,x>})^{-1},x(e+\sqrt{<x,x>})^{-1}>)
$$
$$
= P_\alpha((e+\sqrt{<x,x>})^{-1}<x,x> (e+\sqrt{<x,x>})^{-1})
$$
$$
=P_\alpha(\sqrt{<x,x>} (e+\sqrt{<x,x>})^{-1})^2 \le 1.
$$

b) For each $x \in X$ and each positive number $t$ we have
$$x-x(e+t\sqrt{<x,x>})^{-1}=(x(e+t \sqrt{<x,x>})-x)(e+t\sqrt{<x,x>})^{-1}
$$
$$=tx\sqrt{<x,x>}(e+t\sqrt{<x,x>})^{-1}.$$
Therefore by statement (a) of Lemma 1.2
$$ \bar P_\alpha(x-x(e+t \sqrt{<x,x>})^{-1}) \le t \bar
P_{\alpha}(x\sqrt{<x,x>}).  $$
This implies b). The lemma is proved.

 \begin{cor}
 The set $X^s$ is dense in $X$.
 \end{cor}

{\bf Proof.} Let $x \in X$ be an arbitrary element and $t$
a positive number. Then Lemma 2.3 implies that the elements
$(e+\sqrt{<x,x>})^{-1}, \; x(e+t\sqrt{<x,x>})^{-1}$ belong
to $X^s$ and
$$ \lim \limits_{t \to 0} (x(e+ t\sqrt{<x,x>})^{-1})= x$$
in $X.$ Thus $X^s$ is dense in $X$.

 \begin{th}
$X^s$ is a Hilbert $C^*$-module over the $C^*$-algebra $A^s.$
 \end{th}

{\bf Proof.} First we will show that the restriction of the inner product
$<.,.>$ from $X$ to $X^s$ is an $A^s$-valued inner product on $X^s$.
Indeed, by the Cauchy-Bunyakovskii inequality $(2.1)$ we have
$$
P_\alpha (<x,y>)^2 \le P_\alpha (<x,x>)\cdot P_\alpha (<y,y>)
\le ||<x,x>||^s ||<y,y>||^s
$$
for all $x, y \in X^s$ and $\alpha\in \Delta$.
Therefore, $<x,y> \in A^s$.

Let us prove completeness of $X^s$ with respect to the norm
$$
||x||^s= (||<x,x>)||^s)^\frac12=
\sup\limits_{\alpha \in \bigtriangleup}\bar P_\alpha (x). \eqno(2.4)
$$
Let $\{x_n\}$ be a fundamental sequence in $X^s$, i.e.
for any $\varepsilon> 0$ there exists a natural number $n_\varepsilon$
such that $||x_m-x_n||^s<\varepsilon$ if $m, n> n_\varepsilon$, whence
$$
\bar P_\alpha(x_m-x_n) \le ||x_m-x_n||^s< \varepsilon
$$
for all $\alpha \in \bigtriangleup$ and $m, n> n_\varepsilon$.
This means that $\{x_n\}$ is a Cauchy sequence in $X$ and as
$X$ is complete, so we conclude that the limit
$$ x=\lim\limits_{n \to \infty} x_n $$
exists in $X.$ It follows from the inequalities
$$\left| \bar P_{\alpha}(x_m)-\bar P_{\alpha}(x_n) \right| \le
\bar P_{\alpha}(x_m-x_n), \;
\mid ||x_m||^s-||x_n||^s \mid \le ||x_m-x_n||^s $$
that the sequences $\{\bar P_{\alpha}(x_n)\}$ and
$\{||x_n||^s\}$ are Cauchy sequences of numbers
and therefore are convergent.
Besides, for each $\alpha \in \bigtriangleup$
$$ \bar P_\alpha
(x)=\lim\limits_{n \to \infty} \bar P_\alpha (x_n) \le \lim\limits_{n \to
\infty} ||x_n||^s< \infty$$
hence $x \in X^s$. For all $\alpha \in \bigtriangleup$ and
$n> n_\varepsilon$ we have
$$ \bar P_\alpha (x-x_n)=\lim\limits_{m \to
\infty} \bar P_\alpha (x_m-x_n) \le \varepsilon.$$
Thus for $n>n_\varepsilon$
$$ ||x-x_n||^s=\sup\limits_{\alpha \in \bigtriangleup}\bar
P_\alpha(x-x_n) \le \varepsilon.  $$
This means that the sequence $\{x_n\}$ converges to the element $x$
with respect to the topology of $X^s$. The theorem is proved.

For any $\alpha \in \bigtriangleup$ let us put
$$
I_\alpha=\{a \in A:\;P_\alpha (a)=0\},\qquad J_\alpha=I_\alpha \cap A^s,
$$
$$
\bar I_\alpha=\{x \in X:\;<x,x>\in I_\alpha \},\qquad \bar J_\alpha=
\bar I_\alpha \cap X^s.
$$
The subset $\bar I_\alpha$ is a closed $A$-submodule in $X$,
$J_{\alpha}$ is a closed ideal of the $C^*$-algebra $A^s$ and
$\bar J_\alpha$ is a closed $A^s$-submodule in $X^s$.
Therefore the quotient $X_\alpha:= X/ {\bar I_\alpha}$ is
is a normed space with respect to the norm
$$ ||x+\bar I_\alpha||:=\inf\limits_{y \in \bar
I_\alpha}\bar P_\alpha(x+y) =\bar P_\alpha(x),\qquad x \in X, $$
and the quotient $X^s_\alpha:=X^s/ {\bar J_\alpha}$
is a Banach space with respect to the norm
$$ ||x+\bar J_\alpha||=\inf\limits_{y \in \bar
J_\alpha}||x+y||^s, \qquad x \in X^s.  $$

 \begin{lem}
Let $\{u_\lambda\}_{\lambda \in \Lambda}$ be an approximate
identity in the ideal $J_\alpha$ of the $C^*$-algebra $A^s$.
Then

a)  $\lim\limits_{\lambda}||y-y u_\lambda||^s=0$
for any $y \in \bar J_\alpha.$

b) $||x+\bar J_\alpha||=\lim\limits_{\lambda} ||x-x u_\lambda||^s$
for any $x \in X^s.$

c) $||a+ J_\alpha||=\lim\limits_{\lambda} ||a-a u_\lambda||^s$
for any $a \in A^s.$
 \end{lem}

{\bf Proof.} a) Suppose $y \in \bar J_\alpha$. Then for any
$\beta \in \Delta$ we have
$$
\bar P_\beta(y-y u_\lambda)^2= P_\beta(<y-y u_\lambda,\;y-y u_\lambda>)
=P_\beta((e-u_\lambda)<y,y> (e-u_\lambda))
$$
$$
\le P_\beta(<y,y>-<y,y> u_\lambda)
\le ||<y,y>-<y,y>u_{\lambda}||^s$$
because $P_{\beta}(e-u_{\lambda}) \le ||e-u_{\lambda}||^s \le 1.$
Since $\beta$ is arbitrary, we get
$$(||y-yu_{\lambda}||^s)^2 \le ||<y,y>-<y,y>u_{\lambda}||^s.$$
Using $<y,y> \in J_\alpha$ we obtain $\lim\limits_\lambda
||y- y u_\lambda||^s=0.$

b) Let now $x \in X^s$ and take $\varepsilon> 0$. By definition of
infimum there exists an element $y \in \bar J_\alpha $ such that
$$
||x+y||^s< ||x+\bar J_\alpha||+\frac{\varepsilon}{2}.
$$
Item a) implies that there exists $\lambda_0 \in \Lambda$ such that for
$\lambda> \lambda_0$  one has
$$ ||y-y u_\lambda||^s< \frac{\varepsilon}{2}. $$
Then for all $\lambda> \lambda_0 $
$$ ||x-x u_\lambda||^s=
||x(e-u_\lambda)+y(e-u_\lambda)-y(e-u_\lambda)||^s  \le
||(x+y)(e-u_\lambda)||^s+ ||y-y u_\lambda||^s $$ $$ \le ||(x+y)||^s+||y-y
u_\lambda||^s  \le ||x+\bar
J_\alpha||+\frac{\varepsilon}{2}+\frac{\varepsilon}{2}= ||x+\bar
J_\alpha||+\varepsilon.  $$
Therefore for $\lambda> \lambda_0$ we have
$$ \left| ||x+\bar J_\alpha||-||x-x u_\lambda||^s \right|< \varepsilon. $$
Thus statement b) is true. Statement c) follows from b) for $X=A.$
The lemma is proved.

 \begin{lem}
For any $x \in X^s$ the equality
$$
||x+\bar J_\alpha||= \bar P_\alpha (x)
$$
holds.
 \end{lem}

{\bf Proof.} Let $x \in X^s$ and $\{u_{\lambda}\}_{\lambda \in
\Lambda}$ approximate identity of the ideal $J_{\alpha}.$
Then by statements b) and c) of Lemma 2.4 we have
 $$ ||x+\bar J_{\alpha}||^2=\lim \limits_{\lambda}(||x-x u_{\lambda}||^s)^2=
 \lim \limits_{\lambda}\sup \limits_{\beta} \bar P_\beta(x-x u_{\lambda})^2
 =\lim \limits_{\lambda}\sup\limits_{\beta} P_{\beta}(<x-x u_{\lambda},
 x-xu_{\lambda}>) $$
$$ =\lim\limits_{\lambda} \sup\limits_{\beta}
 P_\beta((e-u_\lambda)<x,x>(e-u_\lambda))
\le \lim\limits_\lambda \sup\limits_\beta P_\beta
(<x,x> -<x,x> u_\lambda)
$$
$$
=\lim\limits_\lambda ||<x,x> -<x,x> u_\lambda||^s=||<x,x>+J_\alpha||
=P_\alpha(<x,x>)=(\bar P_\alpha(x))^2.
$$
The next to last equality follows from uniqueness of $C^*$-norm of
$C^*$-algebra $A^s/J_\alpha$ (see \cite{22}).
Thus we have proved that
$$
||x+\bar J_\alpha|| \le \bar P_\alpha(x).
$$
The inverse inequality is verified immediately:
$$ \bar P_\alpha(x)=||x+\bar I_\alpha||=\inf\limits_{y \in \bar
I_\alpha} \bar P_\alpha(x+y) \le \inf\limits_{y \in \bar J_\alpha} \bar
P_\alpha(x+y)  \le \inf\limits_{y \in \bar J_\alpha}
||x+y||^s=||x+\bar J_\alpha||. $$
The lemma is proved.

 \begin{th}
 The quotient module $X_\alpha=X/\bar I_\alpha$ is a Hilbert
  module over the $C^*$-algebra $A_\alpha$.
 \end{th}

{\bf Proof.}  We define an action of the algebra $A_\alpha$ on $X_\alpha$
by the formula
$$ (x+\bar I_\alpha)(a+I_\alpha):=xa+\bar
I_\alpha, \quad x \in X,\quad a \in A.  $$
With respect to this action $X_\alpha$ is a right $A_\alpha$-module.
Let us define the $A_\alpha$-valued inner product in $X_\alpha$
by the following formula
$$ <x+\bar I_\alpha,\;y+\bar
I_\alpha>:=<x,y>+I_\alpha,\quad x,y \in X. $$
The inner product axioms are easily verified. Positive definiteness
follows from equality (1.1).
Further,
$$ ||<x+\bar I_\alpha,\;x+\bar
I_\alpha>||=||<x,x>+I_\alpha||=P_\alpha(<x,x>)= (\bar
P_\alpha(x))^2=||x+\bar I_\alpha||^2,$$
i.e.
$$ ||x+\bar
I_\alpha||=||<x+\bar I_\alpha,\;x+\bar I_\alpha>||^{\frac12}.$$
To establish the completeness of $X_\alpha$ let us consider the map
$$ \varphi:X^s / \bar J_\alpha \to X/\bar I_\alpha=X_\alpha $$
defined by the formula
$$ \varphi(x+\bar J_\alpha)=x+\bar I_\alpha,\quad x \in X^s.  $$
It is clear that $\varphi$ is an injective linear map.
Besides, it follows from Lemma 2.5 that this map preserves the norm,
  $$
||\varphi(x+\bar J_\alpha)||=||x+\bar I_\alpha||=\bar P_\alpha(x)=
||x+\bar J_\alpha||,\quad x \in X^s.  $$
Consequently, the image $\varphi(X^s /\bar J_\alpha)=\tilde X^s$
is closed in $X_\alpha$. Let us prove that the set $\tilde X^s$
is dense in $X_\alpha$.

Let $x \in X$ be an arbitrary element. Then Lemma 2.3 implies that
for any $t>0\;$ the element $x(e+t\sqrt{<x,x>})^{-1}$ belongs to $X^s$
and
$$
\lim\limits_{t \to 0} x(e+t\sqrt{<x,x>})^{-1}=x
$$
in $X$. Therefore,
$$
x(e+t\sqrt{<x,x>})^{-1}+\bar I_\alpha \in \tilde X^s
$$
and
$$
\lim\limits_{t \to 0} [x(e+t\sqrt{<x,x>})^{-1}+\bar I_\alpha]
=x+\bar  I_\alpha
$$
in $X_\alpha$. Since $\tilde X^s$ is closed in $X_\alpha$,
we conclude that $\tilde X^s=X_\alpha$. Thus the space $X_\alpha$
is complete. The theorem is proved.

 \section{Operators on Hilbert modules over locally $C^*$-algebras}

Let $X$ be a Hilbert module over a locally $C^*$-algebra $A$ and
 $\{ P_\alpha\}_{\alpha\in\Delta}$ a family of $C^*$-seminorms on $A$.

{\bf Definition 3.1.} ${\bf C}$-linear $A$-homomorphism $T:X \to X$
is called {\em bounded operator\/} on the module $X$ if for each
  $\alpha \in \Delta$ there exists a constant $C_\alpha>0$ such that
$$
\bar P_\alpha(Tx) \le C_\alpha \bar P_\alpha(x)
$$
for all $x \in X$.

We denote by ${\rm End}_A(X)$ the set of all bounded operators on $X$.
The following lemma is easily proved.

 \begin{lem}
 For any $\alpha \in \Delta$ the function
$$
\hat P_\alpha(T):=\sup\limits_{\bar P_\alpha(x)\le 1}\bar P_\alpha(Tx)
$$
is a seminorm in ${\rm End}_A(X)$.
 \end{lem}

 \begin{th}
 The set ${\rm End}_A(X)$ is a complte $LMC$-algebra with respect to the
family of seminorms $\{\hat P_\alpha\}_{\alpha \in \Delta}$.
 \end{th}

{\bf Proof.} Let $T_1, T_2 \in {\rm End}_A(X).$ Then for all $\alpha
\in \Delta$ and $x \in X$
$$\bar P_{\alpha}((T_1+T_2)x) \le \bar
P_{\alpha}(T_1x)+\bar P_{\alpha}(T_2x) \le (\hat P_{\alpha}(T_1)+\hat
P_{\alpha}(T_2))\bar P_{\alpha}(x),$$
$$ \hat P_{\alpha}((T_1 T_2)x) \le \hat
P_{\alpha}(T_1)\bar P_{\alpha}(T_2 x) \le \hat P_\alpha(T_1) \hat
P_\alpha(T_2)\bar P_\alpha(x).$$
This implies that $T_1+T_2, \; T_1 T_2 \in
{\rm End}_A(X)$  and
$$\hat P_{\alpha}(T_1+T_2) \le \hat P_{\alpha}(T_1)+\hat
P_{\alpha}(T_2), \qquad \hat P_{\alpha}(T_1T_2) \le \hat P_{\alpha}(T_1)\hat
P_{\alpha}(T_2). \eqno (3.1)$$
Thus ${\rm End}_A(X)$ is an algebra and the function $\hat P_{\alpha}$
is a semi-multiplicative seminorm for this algebra.

Let us verify now
that if $\hat P_\alpha(T)=0$ for all $\alpha\in \Delta$, then $T=0$.

Suppose $\hat P_\alpha(T)=0$ for all $\alpha \in \Delta$.
Then for all $x \in X$ and $\alpha \in \Delta$ we have
$$
\bar P_\alpha (Tx)\le \hat P_\alpha(T)\bar P_\alpha (x)=0,
$$
cosequently, for all $x \in X, \; \alpha \in \Delta$
$$
\bar P_\alpha (Tx)=0,
$$
hence $Tx=0$ for all $x \in X$. This yields $T=0.$

Let us prove completeness of ${\rm End}_A(X)$. Let
$\{T_i\}_{i \in I}\subset {\rm End}_A(X)$ be a Cauchy net,
i.e. for any $\alpha \in \Delta$ and $\varepsilon>0$ there exists $\beta
\in I$ such that for $i_1,i_2>\beta$
$$ \hat P_\alpha(T_{i_1}-T_{i_2})<\varepsilon.\eqno(3.2) $$
Therefore, if $\bar P_\alpha(x) \ne 0$ then for each $\varepsilon>0$
there exists $\beta \in I$ such that for $i_1,i_2>\beta$
$$ \hat P_\alpha(T_{i_1}-T_{i_2})< \frac{\varepsilon}{\bar
P_\alpha(x)}.  $$
Then for $i_1,i_2>\beta$
$$
\bar P_\alpha(T_{i_1}x-T_{i_2}x)\le \bar P_\alpha(T_{i_1}-T_{i_2}) \bar
P_\alpha(x)< \varepsilon. \eqno(3.3)
$$
If $\bar P_\alpha(x)=0$ then $$ \hat P_\alpha(T_{i_1}x-T_{i_2}x)\le \hat
P_\alpha(T_{i_1}-T_{i_2}) \bar P_\alpha(x)=0, $$
i.e. inequality (3.3) holds too.

This yields that for each $x \in X$ \
$\{T_i x\}_{i \in I} \subset X$ is a Cauchy net. Therefore
there exists the limit
$$
T x=\lim\limits_i T_i x.    \eqno(3.4)
$$
It is easy to verify that the map $T$ defined by $(3.4)$ is an
$A$-linear homomorphism on $X$. Let us prove boundedness of $T$.

Since for each $\alpha \in \Delta$ and $i_1,i_2 \in I$
$$
|\hat P_\alpha(T_{i_1})-\hat P_\alpha(T_{i_2})|\le \hat
P_\alpha(T_{i_1}-T_{i_2}), $$
we see that the net $\{\hat P_\alpha(T_i)\}_{i \in I}$
of positive numbers is a Cauchy net. Hence for each $\alpha \in \Delta$
 the limit
$$ \lim\limits_i \hat P_\alpha(T_i)=\lambda_\alpha \ge 0$$
exists.  For all $x \in X,\; \alpha \in \Delta$ we have
$$
\bar P_\alpha(Tx)=\lim\limits_i \bar P_\alpha(T_i x)
\le \lim\limits_i (\hat P_\alpha(T_i)\bar P_\alpha(x))=\lambda_\alpha
\bar P_\alpha(x).
$$
This yields boundedness of $T$. Thus $T \in {\rm End}_A(X)$.

It follows from (3.2) that for $P_\alpha(x) \le 1$ and $i_1,i_2> \beta$
$$ \bar P_\alpha(T_{i_1}x-T_{i_2}x)\le \sup\limits_{\bar
P_\alpha(y) \le 1} \bar P_\alpha((T_{i_1}-T_{i_2})y)
=\hat P_\alpha(T_{i_1}-T_{i_2})< \varepsilon.  $$
Taking limit over $i_2$ we obtain
$$ \bar P_\alpha(T_{i_1} x-Tx)<\varepsilon $$
if $\bar P_\alpha(x) \le 1$ and $i_1>\beta$.
Consequently, for $i_1>\beta$ we have
$$ \hat P_\alpha(T_{i_1}-T)=\sup\limits_{\bar P_\alpha(x) \le 1} \bar
P_\alpha(T_{i_1}x-Tx) \le \varepsilon, $$
whence
$$ \lim T_i= T, $$
i.e. ${\rm End}_A(X)$ is a complete space.
The theorem is proved.

{\bf Definition 3.2.} An operator $T \in {\rm End}_A(X)$ is called
{\em adjointable\/} if there exists an operator $T^* \in {\rm End}_A(X)$
such that
$$<Tx,y>=<x,T^* y>$$
for all $x,y \in X.$

We will denote by ${\rm End}^*_A(X)$ the set of all bounded adjointable
operators on $X$.

 \begin{lem}
  Let $T:X \to X,\; T^*:X \to X$ be maps such that the equality
$$
<x,Ty>=<T^* x,y>
$$
holds for any elements $x$ and $y$ from $X$. Then $T$ and $T^*$ are
bounded operators on $X$, consequently, $T \in {\rm End}_A^*(X)$.
 \end{lem}

{\bf Proof}. For all $x,y,z \in X,\;\lambda \in {\bf C}$ and
$a \in A$ the following equalities hold
$$
<z,T(x+y)>=<T^*z,x+y>=<T^*z,x>+<T^*z,y>
$$
$$
=<z,Tx>+<z,Ty>=<z,Tx+Ty>,
$$
$$
<z,T(\lambda x)>=<T^*z,\lambda x>=<T^*z,x>\lambda=
<z,Tx>\lambda=<z,\lambda Tx>,
$$
$$
<z,T(xa)>=<T^*z,xa>=<T^*z,x>a=<z,Tx>a=<z,(Tx)a>.
$$
As $z$ was arbitrary, so we obtain
$$
T(x+y)=Tx+Ty,\quad T(\lambda x)=\lambda Tx,
$$
$$
T(xa)=(Tx)a.
$$
Thus the desired properties of linearity hold.

Now let us prove that  the submodule $\bar I_\alpha \subset X$
 is invariant with respect to $T$ and $T^*$ for any $\alpha \in \Delta$.
Indeed, if $x \in \bar I_\alpha$, i.e. $\bar P_\alpha(x)=0$,
 then by the Cauchy-Bunyakovskii inequality we have
$$
\bar P_\alpha(Tx)^2=P_\alpha(<Tx,Tx>)= P_\alpha (<x,T^* Tx>)
$$
$$
\le P_\alpha (<x,x>)^\frac12 P_\alpha(<T^*Tx,T^*Tx>)^{\frac12}=
\bar P_\alpha (x)\bar P_\alpha (T^*Tx)=0.
$$
Consequently, $\bar P_\alpha(Tx)=0$, i.e., $Tx \in \bar I_\alpha $.
By the same reason we get $T^* x \in \bar I_\alpha$.

 Thus maps
$$
T_\alpha:X_\alpha \to X_\alpha,  \qquad
T^*_\alpha :X_\alpha  \to X_\alpha
$$
given by the formulas
$$
T_\alpha(x+\bar P_\alpha)=Tx+\bar I_\alpha,\quad
T^*_\alpha(x+\bar I_\alpha)=T^* x+\bar I_\alpha
$$
are well-defined.
 Furthemore,
  $$
<x+\bar I_\alpha,T_\alpha(y+\bar I_\alpha)>=<T^*_\alpha(x+\bar I_\alpha),
y+\bar I_\alpha>
$$
for all $x,y \in X.$
 Since $X_\alpha$ is a Hilbert module over the $C^*$-algebra $A_\alpha$,
 we obtain that $T_\alpha$ is a bounded $A_\alpha$-homomorphism
on $X_\alpha$ (see \cite{10,17}), i.e. there exists a constant $C_\alpha>0$
such that
  $$ ||T_\alpha (x+\bar I_\alpha)|| \le C_\alpha
||x+\bar I_\alpha||\eqno(3.5) $$
for all $x \in X$. It follows from (3.5) that
 $$ \bar P_\alpha(Tx) \le C_\alpha \bar P_\alpha(x)$$
for all $x \in X$ since $\bar P_\alpha(Tx)=||Tx+\bar I_\alpha||=
||T_\alpha(x+\bar I_\alpha)||$
 and
 $||x+\bar I_{\alpha}||=\bar P_{\alpha}(x),$
 whence, since $\alpha \in \Delta$ is arbitrary, it follows that
 the operator $T$ is bounded. The lemma is proved.

  \begin{th}
  The set ${\rm End}^*_A(X)$ is a locally $C^*$-algebra with respect to the
family of seminorms $\{\hat P_\alpha\}_{\alpha \in \Delta}$.
  \end{th}

{\bf Proof.} It is easy to verify that if $T,Q \in {\rm End}^*_A(X)$ and
$\lambda \in C$, then

$$(T+Q)^*=T^*+Q^*, \quad (\lambda T)^*=\bar \lambda T^*,$$

$$(TQ)^*=Q^*T^*,\quad (T^*)^*=T.$$
Thus ${\rm End}^*_A(X)$ is a $*$-algebra.
By statement (3) of Lemma 2.2 for each $T \in {\rm End}^*_A(X)$ we have
$$
\hat P_\alpha(T^* T)=\sup\limits_{\bar P_\alpha(x)\le 1}\bar P_\alpha(T^*
 Tx)= \sup\limits_{\bar P_\alpha(x)\le 1} \sup\limits_{\bar P_\alpha(y)\le
1}P_\alpha(<T^*Tx,y>) $$ $$ =\sup\limits_{\bar P_\alpha(x)\le 1}
\sup\limits_{\bar P_\alpha(y)\le 1} P_\alpha(<Tx,Ty>) \ge
\sup\limits_{\bar P_\alpha(x) \le 1}P_\alpha(<Tx,Tx>) 
=\sup\limits_{\bar P_\alpha(x)\le 1}\bar P_\alpha(Tx)^2=\hat P_\alpha(T)^2,
$$
i.e.
$$
\hat P_\alpha(T^* T)\ge \hat P_\alpha(T)^2.  \eqno(3.6)
$$
It follows from (3.1) and (3.6) that
$$
\hat P_\alpha(T^* T)=\hat P_\alpha(T)^2,
$$
whence one can easy obtain that
$$
\hat P_\alpha(T^*)=\hat P_\alpha(T).
$$
Thus $\hat P_{\alpha}$ is a $C^*$-seminorm on ${\rm End}_A^*(X)$
for any $\alpha \in \Delta$.

It remains to prove completeness of ${\rm End}^*_A(X)$.
 For this aim it is enough  to verify that ${\rm End}^*_A(X)$ is closed in
${\rm End}_A(X)$. Let $\{T_i\}_{i \in I} \subset {\rm End}^*_A(X)$ and
$\lim T_i=T.$ Then $\{T^*_i\}_{i \in I}$ is a Cauchy net in
${\rm End}_A(X)$ and therefore converges to some operator $Q \in {\rm
End}_A(X)$.
For all $x,y \in X$ and $i \in I$ we have
$$ <T_i x,y>=<x,T^*_i y> $$
therefore taking limit over $i$ we obtain
$$ <Tx,y>=<x,Qy>.  $$
Hence the operator $T$ is adjointable and $T^*=Q$.  Consequently,
$T \in {\rm End}^*_A(X)$. The theorem is proved.

For any $\alpha \in \Delta$ and $T \in {\rm End}^*_A(X)$ we define
the map $T_\alpha:X_\alpha\to X_\alpha$ by the formula
$$ T_\alpha(x+\bar I_\alpha):=Tx+\bar I_\alpha.$$
One can easy verify that $T_\alpha \in {\rm End}^*_{A_\alpha}(X_\alpha)$
(see the proof of Lemma 3.2) and the map
$$ T \mapsto T_\alpha $$
is a $*$-homomorphism from the $LC^*$-algebra ${\rm End}^*_A(X)$ to
the $C^*$-algebra ${\rm End}^*_{A_\alpha}(X_\alpha)$. Besides,
$$ ||T_\alpha||=\hat P_\alpha(T).  $$

 \begin{prop}
  Let $T \in {\rm End}^*_A(X)$. Then

$a)\;Sp(T)=\bigcup_{\alpha \in \Delta} Sp(T_\alpha);$

$b)\; T$ is self-adjoint if and only if $T_\alpha$
is self-adjoint for any $\alpha \in \Delta;$

$c)\; T$ is a positive element of the $LC^*$-algebra ${\rm End}^*_A(X)$
if and only if $T_\alpha$ is a positive element of the
$C^*$-algebra ${\rm End}^*_{A_\alpha}(X_\alpha)$ for any  $\alpha \in
\Delta$.
 \end{prop}

 \begin{prop}
  For an operator $T:X \to X$ the following conditions are equivalent:

1) $T$ is a positive element of the $LC^*$-algebra ${\rm End}^*_A(X)$;

2) for any element $x \in X$ the inequality $<Tx,x> \ge 0$ holds,
i.e. this element is positive in $A$.
\end{prop}

One can prove Proposition 3.2 similarly to Lemma 4.1 of \cite{8}
or by using statement c) of Proposition 3.1.
\begin{prop}
  The map $T \mapsto T|_{X^s}$ is an isomorphism of the $C^*$-algebras
$({\rm End}^*_A(X))^s$ and ${\rm End}^*_{A^s}(X^s)$.
 \end{prop}

{\bf Acknowledgement.}
Authors are grateful to A.~S.~Mishchenko, E.~V.~Troitsky,
V.~M.~Manui\-lov and A.~A.~Pavlov for useful discussions and remarks.
We are also grateful to R.~Meyer and N.~C.~Phillips for pointing out the
papers \cite{a,a1,b}.

\end{document}